\documentclass[letterpaper, 10 pt, journal, twoside]{ieeetran}
\IEEEoverridecommandlockouts
\usepackage{theorem}%
\usepackage{amsmath,amssymb,amsfonts}%
\usepackage{graphicx}      
\usepackage{color,dsfont}

\usepackage{enumitem} 
\definecolor{mycol}{rgb}{0,0,1}
\definecolor{mcc}{rgb}{0,0.4,0.6}
\usepackage{amsmath,amssymb}
\usepackage[mathscr]{eucal}
\usepackage{times,pifont}
\usepackage{framed}
\usepackage[bookmarks,bookmarksopen,bookmarksdepth=3]{hyperref}
\usepackage{xcolor}

\usepackage{tikz}
\usetikzlibrary{arrows}
\usetikzlibrary{graphs,decorations.pathmorphing,decorations.markings}
\usetikzlibrary{calc}

\theoremstyle{change}%
\newtheorem{definition}{Definition:}[section]%
\newtheorem{theorem}[definition]{Theorem}%
\newtheorem{lemma}[definition]{Lemma}%
\newtheorem{corollary}[definition]{Corollary}%

{\theorembodyfont{\rmfamily}}%
{\theorembodyfont{\rmfamily}\newtheorem{example}[definition]{Example}}%

\newcommand{\N}{\mathbb{N}}%
\newcommand{\R}{\mathbb{R}}%
\newcommand{\cQ}{\mathcal{Q}}%
\newcommand{\tm}{\times}%
\DeclareMathOperator{\rank}{rank}

\DeclareMathOperator{\range}{ran}

\allowdisplaybreaks%

\title{\Huge A case study of port-Hamiltonian
systems with a moving interface}

\author{
Alexander Kilian\thanks{A.~Kilian, A.~Mironchenko, F.~Wirth are with the Faculty of Computer Science and Mathematics, University of Passau, 94032 Passau, Germany; e-mail: \texttt{$\{$alexander.kilian,andrii.mironchenko,fabian.} \texttt{(lastname)$\}$@uni-passau.de}. A.~Kilian is supported by BMBF through the grant 16ME0619. A.~Mironchenko is supported by the DFG through the grant MI 1886/2-2.},
Bernhard Maschke\thanks{B.~Maschke is with LAGEPP, UCB Lyon 1 - CNRS UMR 5007, CPE Lyon - Bâtiment 308 G, Université Claude Bernard Lyon-1, 43, bd du 11 Novembre 1918, F-69622 Villeurbanne cedex, France; e-mail: \texttt{bernhard.maschke@univ-lyon1.fr}.},
Andrii~Mironchenko,
Fabian Wirth
}

\begin{document}

\maketitle%
\thispagestyle{empty} 

\begin{abstract}
	We model two systems of two conservation laws defined on complementary spatial intervals and coupled by a moving interface as a single non-autonomous port-Hamiltonian system, and provide sufficient conditions for its Kato-stability. An example shows that these conditions are quite restrictive. The more general question under which conditions an evolution family is generated remains open.
\end{abstract}


%

	\section{Introduction}

The port-Hamiltonian framework \cite{Duindam2009, Jacob2012, Schaft2014} 
is a systematic theory for modeling, analysis and control of \emph{physical systems} belonging to mechanical, electric, hydraulic, thermal, and other domains.
It rests on two eponymous concepts: the system's \emph{Hamiltonian} (that is, its total energy) stemming from mechanics \cite{Arnold1978}, and the \emph{port-based modeling approach} \cite{Duindam2009}, which is based on the power-conserving interconnection of various system components by means of ports. 


Following the semigroup approach \cite{Curtain2020}, port-Hamiltonian distributed parameter systems can be represented as abstract Cauchy problems of the form
\begin{align*}
	\dot{x}(t) &= A_{\mathcal{Q}}x(t), \quad t > 0, \\
	x(0) &= x_0 \in X,
\end{align*}
where $A_{\mathcal{Q}} \colon D(A_{\mathcal{Q}}) \subset X \to X$ is a linear unbounded operator on a Hilbert space $X$ (the so-called \emph{port-Hamiltonian operator}) associated with the port-Hamiltonian system. 
The port-Hamiltonian operator stores information concerning the interdomain coupling of the multi-physics system, as well as the mechanical and physical properties. 
Frequently, $A_{\mathcal{Q}} = \mathcal{J}(\mathcal{Q}\cdot)$, where $\mathcal{J}$ is a formally skew-symmetric matrix differential operator modeling the interdomain coupling, and $\mathcal{Q}$ is a coercive and bounded multiplication operator representing the physical properties  \cite[Section 3]{Gorrec2005}, \cite[pp. 226-228]{Duindam2009}.
As discussed in \cite{KMMW23,Kilian2022}, a significant research effort was devoted to the analysis of port-Hamiltonian systems connected via interfaces.
Such systems appear, e.g., due to the distinct physical properties on the respective subdomains. 
To cover couplings of systems with a non-fixed boundary between them in port-Hamiltonian terms, in \cite{Diagne2013} the authors proposed the port-Hamiltonian formulation for two systems that are coupled by a moving interface. However, \cite{Diagne2013} was concentrated on the modeling, and the well-posedness of the obtained model has not been investigated in \cite{Diagne2013}.

In \cite{KMMW23}, we proposed a mathematically rigorous approach for the study of port-Hamiltonian systems with an interface as abstract evolution equations. Using the framework of \cite{Jacob2012}, easy-to-check criteria for boundary and interface conditions that guarantee that the system with a fixed interface induces a contraction semigroup were given in \cite{KMMW23}. 
Under somewhat stronger conditions, exponential stability of such a system was established in \cite{KMMW23}. These findings have been applied in \cite{KMMW23} to two acoustic waveguides coupled by an interface consisting of some membrane and in \cite{Kilian2022} to the interface coupling of two lossless transmission lines.

\textbf{Contribution.} 
In \cite{KMMW23} only systems with a stationary interface have been considered. 
In this note, we show that the approach advocated in \cite{KMMW23} is 
applicable also to systems with a moving interface. 
 In Section \ref{Subsection Port-Hamiltonian Formulation}, we consider a port-Hamiltonian formulation of the system with an interface. We introduce the formally skew-symmetric differential operator associated with the aggregate system. Furthermore, we define the boundary and interface port variables in case of a fixed interface position, and state an energy balance equation for this class of systems. 
In Section \ref{Section Stability of the Family of Infinitesimal Generators}, we shall define a family of port-Hamiltonian operators that keeps track of the position of the moving interface, and analyze the resulting time-variant evolution problem associated with this family. We present sufficient conditions for the Kato-stability of this family in the sense of Pazy \cite[Section 5.2]{Pazy1983}. This is a key property for the further analysis of the system formulation suggested in \cite{Diagne2013} and \cite{Kilian2022}. 

	\textbf{Notation.}
	Throughout the paper, $X$ is a Hilbert space and $\mathcal{L}(X)$ is the space of linear bounded operators on $X$.
	We denote by $\mathcal{C}([a,b], X)$ and $\mathcal{C}^k([a,b], X)$ the vector spaces of $X$-valued, continuous and $k$-times continuously differentiable functions $f \colon [a,b] \to X$. Also, $L^2([a,b], \mathbb{R}^n)$ denotes the vector space of $\mathbb{R}^n$-valued square-integrable functions on the interval $[a,b]$, and $L^{\infty}([a,b], \mathbb{R}^n)$ is the vector space of $\mathbb{R}^n$-valued essentially bounded functions on $[a,b]$. Moreover, $H^1([a,b], \mathbb{R}^n)$ denotes the Sobolev space of order $1$ on $[a,b]$. 
 For real, square, symmetric matrices $A,B$ the expression $A\leq B$ means that $B-A$ is symmetric positive semi-definite. For functions $y:[a,b]\setminus \{l\} \to \R^n$ the notations $y(l^+), y(l^-)$ denote the right (resp. left) limit of $y$ at $l$, assuming existence.

\section{Port-Hamiltonian Formulation}
\label{Subsection Port-Hamiltonian Formulation}

Let $a<0<b$ and $l \in (a,b)$ be the fixed interface position. 
Consider two systems of two scalar conservation laws each of the form
\begin{align}
	\begin{split}
		\hspace{-3mm}\frac{\partial x^-}{\partial t} (z,t) &= P_1 \frac{\partial}{\partial z} \big(\mathcal{Q}^-(z) x^-(z,t)\big), \ \  z \in [a,l), \; t >0,  \\
		\hspace{-3mm}\frac{\partial x^+}{\partial t} (z,t) &= P_1 \frac{\partial}{\partial z} 
		\big(\mathcal{Q}^+(z) x^+(z,t)\big), \  \ z \in (l,b], \; t > 0,
	\end{split} 
	\label{Systems of Conservation Laws} 
\end{align}
defined on the respective interval $[a,l)$ or $(l,b]$, where
\begin{equation}
    P_1 = \begin{bmatrix}
        0 & -1 \\
        -1 & 0
    \end{bmatrix} \in \mathbb{R}^{2 \times 2},
    \label{Matrix P1}
\end{equation}
and $\mathcal{Q}^{\pm} \in L^{\infty}([a,b], \mathbb{R}^{2 \times 2})$ are pointwise symmetric and satisfy $mI_2 \leq \mathcal{Q}^{\pm}(z) \leq MI_2$ for for suitable $0 < m \leq M$ and a. e. $z \in [a,b]$. Here $I_2$ is the $2\tm 2$ identity matrix.

In the following, we propose a port-Hamiltonian formulation of the model \eqref{Systems of Conservation Laws} on the composed spatial domain $[a,b]$. 
We introduce state variables defined on the composed domain, and we need to define \emph{interface relations} at $z = l$. For a detailed discussion, we refer to \cite{Diagne2013}, \cite[Section 5.2]{Kilian2022}, \cite{KMMW23}. 

	Consider the state space $X = L^2([a,b], \mathbb{R}^2)$ and the system
	\begin{equation}
		\label{Simplified PH-System}
		\frac{\partial}{\partial t} x(z,t) = \mathcal{J}_l \big( \mathcal{Q}_l(z) x(z,t) \big), \quad z \in [a,b] \setminus \{l\}, \; t >0, 
	\end{equation}
	where $\mathcal{J}_l$ is defined below, and 
 \begin{equation}
 \mathcal{Q}_l = c_l^- \mathcal{Q}^- + c_l^+ \mathcal{Q}^+ \in \mathcal{L}(X)
 \label{eq:Q_l}
\end{equation}
 is a coercive matrix multiplication operator. The functions $c_l^-$ and $c_l^+$ are called \emph{color functions} and are used to keep track of the interface position:
 \begin{equation}
		c_l^-(z) = \begin{cases} 
			1, & z \in [a,l), \\
			0, & z \in [l,b],
		\end{cases}  \quad 
            c_l^+(z) =  \begin{cases} 
			0, & z \in [a,l], \\
			1, & z \in (l,b].
		\end{cases}
		\label{Characteristic Functions c- and c+}
	\end{equation}
Since, for now, we assume that the interface position is fixed, the color functions do not depend on time.

	We endow $X$ with the inner product $\langle \cdot, \cdot \rangle_{\mathcal{Q}_l} = \frac{1}{2} \langle \cdot, \mathcal{Q}_l \cdot \rangle_{L^2}$: 
	\begin{equation}
		\langle x , y \rangle_{\mathcal{Q}_l} = \frac{1}{2} \int_{a}^{b} y^{\top}(z) \mathcal{Q}_l(z) x(z) \, dz, \quad x, y \in X,
		\label{Inner Product wrt Q0}
	\end{equation}
	with associated norm $\| \cdot \|_{\mathcal{Q}_l}$.
	Note that $\langle \cdot , \cdot \rangle_{\mathcal{Q}_l}$ is equivalent to the standard inner product $\langle \cdot , \cdot \rangle_{L^2}$ on $X$. 
	The Hamiltonian $H \colon X \to \mathbb{R}$ on the energy space $(X, \| \cdot \|_{\mathcal{Q}_l})$ is defined as
	\begin{equation}
		H(x) = \frac{1}{2} \int_{a}^{b} x^{\top}(z) \mathcal{Q}_l(z) x(z) \, dz =  \|x\|_{\mathcal{Q}_l}^2, \quad x \in X.
		\label{Simplified PH-System - Hamiltonian}
	\end{equation}
	The operator $\mathcal{J}_l \colon D(\mathcal{J}_l) \subset X \to X$ is given by
	\begin{align}
		\hspace*{-0.15 cm}	D(\mathcal{J}_l) &= \left\{ x = (x_1,x_2) \in X \mid x_1 \in D(\mathbf{d}_l^{\ast}), \, x_2 \in D(\mathbf{d}_l)\right\}, \nonumber \\
	\hspace*{-0.15 cm}		\mathcal{J}_l x &= \begin{bmatrix}
				0 & \mathbf{d}_l  \\
				- \mathbf{d}_l^{\ast} & 0 
			\end{bmatrix} \begin{bmatrix}
				x_1 \\
				x_2
			\end{bmatrix}
			=
			 \begin{bmatrix}
				\mathbf{d}_l x_2 \\
				- \mathbf{d}_l^{\ast} x_1
			\end{bmatrix}, \quad x \in D(\mathcal{J}_l),
		\label{Operator Jl}
	\end{align}
	where the operator $\mathbf{d}_l$ is given by
    \begin{align}
    \begin{split}
    	 D(\mathbf{d}_l) & = H^1([a,b], \mathbb{R}),  \\
    	\mathbf{d}_l x &= - \left[ \frac{d}{dz} (c_l^- x) + \frac{d}{dz} (c_l^+ x) \right], \quad x \in D(\mathbf{d}_l),
    \end{split}
    \label{Operator dl}
    \end{align}
    and its formal adjoint (see lines after \eqref{Relation dl and dl_ast}) $\mathbf{d}_l^{\ast}$ is given by
    \begin{align}
    \begin{split}
	D(\mathbf{d}_l^{\ast}) &= \left\{ y \in L^2([a,b], \mathbb{R}) \mid y_{|(a,l)} \in H^1((a,l), \mathbb{R}), \right. \\
 &\qquad \qquad \qquad \qquad \quad \hspace{0.1 cm} \left. y_{|(l,b)} \in H^1((l,b), \mathbb{R}) \right\}, \\
	\mathbf{d}_l^{\ast}y &= \left( - \mathbf{d}_l - \left[ \frac{d}{dz} c_l^- + \frac{d}{dz} c_l^+ \right] \right) y, \quad y \in D(\mathbf{d}_l^{\ast}).
 \end{split}
 \label{Operator dl_ast}
    \end{align}
    For functions $x = c_l^- x^- + c_l^+ x^+ \in D(\mathbf{d}_l)$ and $y = c_l^- y^- + c_l^+ y^+ \in D(\mathbf{d}_l^{\ast})$, we have in particular that
    \begin{align}
    \mathbf{d}_l x &= - c_l^- \frac{d}{dz} x^- - c_l^+ \frac{d}{dz} x^+, \\
    \mathbf{d}_l^{\ast} y &= c_l^- \frac{d}{dz} y^- + c_l^+ \frac{d}{dz} y^+.
    \end{align} 
    For all $x \in D(\mathbf{d}_l), \,  y \in D(\mathbf{d}_l^{\ast})$, we have the following relation between the operators $\mathbf{d}_l$ and $\mathbf{d}_l^{\ast}$:
 \begin{equation}
 \begin{split}
		\quad \langle \mathbf{d}_l x , y \rangle_{L^2} &=  - \big[ x(z) y(z) \big]_{a}^{b}  \\
  &\quad + x(l) \left[ y(l^+) - y(l^-) \right] + \langle x , \mathbf{d}_l^{\ast} y \rangle_{L^2}.
  \end{split}
  \label{Relation dl and dl_ast}
	\end{equation} 
    In particular, equation \eqref{Relation dl and dl_ast} describes the \emph{formal adjointness} of $\mathbf{d}_l$ and $\mathbf{d}_l^{\ast}$, that is, adjointness up to boundary and interface evaluations.

	On any subinterval of $[a,b]$ not containing the interface position $l \in (a,b)$, the operator $\mathcal{J}_l$ simply acts as the matrix differential operator $P_1 \frac{d}{dz}$, where $P_1$ is defined in \eqref{Matrix P1}. Thus, \eqref{Simplified PH-System} is a reformulation of \eqref{Systems of Conservation Laws}. Next, we reveal the port-Hamiltonian structure of \eqref{Simplified PH-System}.
	As port-Hamiltonian models are extensions of Hamiltonian systems, it is required that the associated differential operator is formally skew-symmetric.
	\begin{lemma}[\cite{KMMW23}, Lemma 2.3]
		The operator $\mathcal{J}_l \colon D(\mathcal{J}_l) \subset X \to X$ defined in \eqref{Operator Jl} is formally skew-symmetric satisfying for all $x = (x_1, x_2)$ and $y = (y_1, y_2) \in D(\mathcal{J}_l)$,
  \begin{align}
  \begin{split}
&\hspace*{-0.15 cm}\langle \mathcal{J}_l x , y \rangle_{L^2} + \langle x , \mathcal{J}_l y \rangle_{L^2} = \big[ x^{\top}(z) P_1 y(z) \big]_{a}^{b} \\
&\hspace*{-0.2 cm}+ x_2(l) \left[ y_1(l^+) - y_1(l^-) \right]  + y_2(l) \left[ x_1(l^+) - x_1(l^-) \right].
  \end{split}
  \label{Skew-Symmetry of J0}
  \end{align}
	\end{lemma}

	
	Next, we augment the system \eqref{Simplified PH-System} with boundary and interface port variables. 
For all $x \in X$ with $ \mathcal{Q}_l x \in D(\mathcal{J}_l)$, the \emph{boundary flow} $f_{\partial} = f_{\partial, \mathcal{Q}_lx} \in \mathbb{R}^2$ and the \emph{boundary effort} $e_{\partial} = e_{\partial, \mathcal{Q}_l x} \in \mathbb{R}^2$ is 
	\begin{equation}
		\label{Boundary Flow and Effort}
		\begin{bmatrix}
			f_{\partial, \mathcal{Q}_l x}\\
			e_{\partial, \mathcal{Q}_l x}
		\end{bmatrix} = \begin{bmatrix}  \frac{1}{\sqrt{2}} P_1 [(\mathcal{Q}_l x)(b) -(\mathcal{Q}_lx)(a)] \\
			\frac{1}{\sqrt{2}} [(\mathcal{Q}_l x)(b) + (\mathcal{Q}_lx)(a))]
		\end{bmatrix},
	\end{equation} 
and the \emph{interface flow} $f_I = f_{I, \mathcal{Q}_l x} \in \mathbb{R}$ and the \emph{interface effort} $e_{I} = e_{I, \mathcal{Q}_l x} \in \R$ are given as follows:
	\begin{align}
		f_{I, \mathcal{Q}_l x} &= \left( \mathcal{Q}_l x \right)_2 (l^+) = \left( \mathcal{Q}_l x \right)_2(l^-), \label{Continuity Equation} \\
		-e_{I, \mathcal{Q}_l x} &= \left( \mathcal{Q}_l x \right)_1(l^+) - \left( \mathcal{Q}_l x \right)_1 (l^-). \label{Balance Equation}
	\end{align}
	For all $ \mathcal{Q}_lx, \mathcal{Q}_ly \in D(\mathcal{J}_l)$, we may write the last expression of \eqref{Skew-Symmetry of J0} with respect to the boundary and interface port variables as follows: 
	\begin{align}
		\begin{split}
			&\quad \big[ (\mathcal{Q}_lx)^{\top}(z) P_1 (\mathcal{Q}_ly)(z) \big]_{a}^{b} \\
   &\quad + \left( \mathcal{Q}_l x \right)_2(l) \left[ \left( \mathcal{Q}_l y \right)_1(l^+) - \left( \mathcal{Q}_l y \right)_1(l^-) \right]   \\
			&\quad + \left( \mathcal{Q}_l y \right)_2(l) \left[ \left( \mathcal{Q}_l x \right)_1(l^+) - \left( \mathcal{Q}_l x \right)_1(l^-) \right] \\
			&\qquad\qquad= \langle e_{\partial, \mathcal{Q}_l y} , f_{\partial, \mathcal{Q}_l x} \rangle + \langle e_{\partial, \mathcal{Q}_l x} , f_{\partial, \mathcal{Q}_l y} \rangle \\
   &\qquad\qquad\qquad\qquad - f_{I, \mathcal{Q}_l x} e_{I, \mathcal{Q}_l y} - f_{I, \mathcal{Q}_l y} e_{I, \mathcal{Q}_l x}. 
		\end{split}
		\label{Auxiliary Equation Power Pairing}
	\end{align}
	Let us substitute the relation \eqref{Auxiliary Equation Power Pairing} into equation \eqref{Skew-Symmetry of J0}. This yields for all $\mathcal{Q}_l x, \mathcal{Q}_l y \in D(\mathcal{J}_l)$, 
\begin{align}
		&\langle \mathcal{J}_l (\mathcal{Q}_l x), \mathcal{Q}_l y \rangle_{L^2} + \langle \mathcal{Q}_l x, \mathcal{J}_l (\mathcal{Q}_l y) \rangle_{L^2} =\langle e_{\partial, \mathcal{Q}_l y} , f_{\partial, \mathcal{Q}_l x} \rangle \nonumber \\
  &+ \langle e_{\partial, \mathcal{Q}_l x} , f_{\partial, \mathcal{Q}_l y} \rangle - f_{I, \mathcal{Q}_l x} e_{I, \mathcal{Q}_l y} - f_{I, \mathcal{Q}_l y} e_{I, \mathcal{Q}_l x}.
		\label{Skew-Symmetry of J0 - Representation 2}
	\end{align}
	
	Now we express a (power) balance equation for classical solutions $x \in C^1([0, \infty), X)$ of the first-order system \eqref{Simplified PH-System} in terms of the boundary and interface port variables.
	\begin{lemma}
		\label{Lemma Balance Equation}
		Let $x \in C^1([0, \infty), X)$ be a classical solution of the system \eqref{Simplified PH-System} with Hamiltonian \eqref{Simplified PH-System - Hamiltonian}. Then for all $t \geq 0$, 
		\begin{equation}
			\label{Power Balance Equation wrt Boundary and Interface Port}
			\hspace*{-2mm}\frac{d}{dt} \|x(t)\|_{\mathcal{Q}_l}^2\! =\!	\frac{d}{dt}H(x(t))\! =\!  \langle e_{\partial}(t), f_{\partial}(t) \rangle - e_I(t) f_I(t),
		\end{equation}
where 
\[
e_{\partial}(t) = e_{\partial, \mathcal{Q}_lx(t)}, \quad f_{\partial}(t) = f_{\partial, \mathcal{Q}_lx(t)},
\]
\[
e_I(t) = e_{I, \mathcal{Q}_lx(t)}, \quad f_I(t) = f_{I, \mathcal{Q}_lx(t)}.
\]
 \end{lemma}
For a proof see \cite[Lemma 2.4]{KMMW23}.

	Equation \eqref{Power Balance Equation wrt Boundary and Interface Port} states that the change of energy is equal to the power flow both at the boundary and at the interface position.

In \cite[Section 5.4]{Kilian2022}, \cite{KMMW23}, the system \eqref{Simplified PH-System} has been studied by considering the corresponding abstract Cauchy problem 
\begin{align}
	\begin{split}
		\dot{x}(t) &= A_{\mathcal{Q}_l}x(t),  \quad t > 0, \\
			x(0) &= x_0 \in D(A_{\mathcal{Q}_l}),
	\end{split}
			\label{Cauchy Problem}
\end{align}
with suitable boundary and interface conditions defined with respect to the boundary and interface port variables, respectively. The port-Hamiltonian operator $A_{\mathcal{Q}_l} \colon D(A_{\mathcal{Q}_l}) \subset X \to X$ associated with \eqref{Simplified PH-System} is defined as follows:
\begin{align}
			D(A_{\mathcal{Q}_l}) &= \left\{ x \in X \, \big| \, \mathcal{Q}_lx \in D(\mathcal{J}_l), \, f_{I, \mathcal{Q}_lx} = r e_{I, \mathcal{Q}_lx}, \right. \nonumber \\
   &\qquad \qquad \qquad \qquad \qquad \quad \left. W_B \begin{bmatrix}
				f_{\partial, \mathcal{Q}_l x} \\
				e_{\partial, \mathcal{Q}_l x}
			\end{bmatrix} = 0 \right\}, \nonumber \\
			A_{\mathcal{Q}_l}x &= \mathcal{J}_l (\mathcal{Q}_lx), \quad  x \in D(A_{\mathcal{Q}_l}),
		\label{def:A_Q}
	\end{align}
where $W_B \in \R^{2 \times 4}$ and $r \in \R$. It has been characterized under which boundary and interface conditions (i.e., for which $W_B$ and $r$) the operator $A_{\mathcal{Q}_l}$ generates a contraction semigroup (cf. \cite[Theorem 5.14]{Kilian2022}, \cite[Theorem 4.4]{KMMW23}). In particular, this guarantees well-posedness of the system \eqref{Cauchy Problem} in the sense of \cite[Definition II.6.8]{Engel2000}. Furthermore, sufficient conditions for the exponential stability of \eqref{Cauchy Problem} have been presented (cf. \cite[Section 5.4.3]{Kilian2022}, \cite[Section 5]{KMMW23}). These results hold for every stationary interface position $l \in (a,b)$. 

In the next section, we make the first steps towards analyzing a more complex case of port-Hamiltonian systems with a moving interface.

	\section{Stability of the Family of Infinitesimal Generators}
	\label{Section Stability of the Family of Infinitesimal Generators}
	In this section, we shall analyze the port-Hamiltonian system \eqref{Simplified PH-System} in case of a moving interface. We define a family of port-Hamiltonian operators that encompasses the position of the interface, and associate this family with an evolution problem of the form
	\begin{align}
	\begin{split}
		\dot{x}(t) &= A(t)x(t), \quad 0 \leq s < t \leq \tau, \\
		x(s) &= x_0 \in X.
	\end{split}
	\label{Evolution Problem}
	\end{align}
	We shall impose conditions guaranteeing that this family is stable in the following sense.	
	\begin{definition}[Kato-stability; {\cite[Chapter 5]{Pazy1983}}]
		\label{Definition Stable Family of Generators}
		A family $(A(t))_{t \in [0, \tau]}$ of infinitesimal generators of $C_0$-semigroups on $X$ is called \emph{Kato-stable} (or just \emph{stable}) if there exist constants $M \geq 1$ and $\omega \in \mathbb{R}$, called stability constants, such that the following holds:
		\begin{enumerate}
			\item[(i)] For $0 \leq t \leq \tau$ we have $(\omega, \infty) \subset \rho(A(t))$.
			\item[(ii)] For every finite sequence $0 \leq t_1 \leq t_2 \leq \ldots \leq t_k \leq \tau$, $k \in \mathbb{N}$, it holds that
			\begin{equation}
				\label{Definition Stable Family of Generators - Condition 2}
				\left\| \prod_{j = 1}^{k} R(\lambda, A(t_j)) \right\|_{\mathcal{L}(X)} \leq \frac{M}{(\lambda - \omega)^k}, \quad \lambda > \omega.
			\end{equation}
		\end{enumerate}
	\end{definition}
Here $\rho(A)$ denotes the resolvent set of an operator $A$, and for $\lambda \in\rho(A)$, $R(\lambda, A)$ denotes the resolvent of $A$ at $\lambda$. 
 
	Recall that the model \eqref{Simplified PH-System} can be defined for arbitrary, but fixed interface positions $l \in (a,b)$ by changing the associated operators accordingly. Unfortunately, the state space $(X, \|\cdot\|_{\mathcal{Q}_l})$ has to be changed for different interface positions $l \in (a,b)$ to guarantee that the aforementioned statements hold true. However, we have the following norm equivalence between $\| \cdot \|_{\mathcal{Q}_0}$ and $\| \cdot \|_{\mathcal{Q}_l}$ for $l \in (a,b)$: for all $x \in X$ we have
	\begin{equation}
		\label{Norm Equivalence}
		\frac{m}{M} \| x \|_{\mathcal{Q}_0}^2 
  \leq \| x\|_{\mathcal{Q}_l}^2 
  \leq \frac{M}{m} \|x \|_{\mathcal{Q}_0}^2.
	\end{equation}
	Note that the constants $\sqrt{\frac{m}{M}}$ and $\sqrt{\frac{M}{m}}$ are independent of the interface position $l$.
	
 Let  $0 < \tau < \infty$, and let $l \colon [0, \tau] \to (a,b)$ be continuously differentiable. Let $(A(t))_{t \in [0, \tau]}$ be a family of port-Hamiltonian operators, where  $A(t) := A_{\mathcal{Q}_{l(t)}}\colon D(A(t))\subset X \to X$ is given by
	\begin{align}
        \begin{split}
			\hspace*{-0.2 cm}D(A(t)) = \Big\{ x \in X \mid & \mathcal{Q}_{l(t)} x \in D(\mathcal{J}_{l(t)}),\\
	\hspace*{-0.2cm}		& f_I(t) = re_I(t), \, W_B \begin{bmatrix}
				f_{\partial}\\
				e_{\partial}
			\end{bmatrix} = 0 \Big\},
        \end{split}
		\label{Operator Family A(t)-I}
	\end{align} 
where 
\begin{equation*}
\begin{aligned}
e_{\partial} = e_{\partial}(t) &= e_{\partial, \mathcal{Q}_{l(t)}x}, \\
e_I(t) &= e_{I, \mathcal{Q}_{l(t)}x},
\end{aligned}
\qquad
\begin{aligned}
f_{\partial} = f_{\partial}(t) &= f_{\partial, \mathcal{Q}_{l(t)}x}, \\
f_I(t) &= f_{I, \mathcal{Q}_{l(t)}x}
\end{aligned}
\end{equation*}
are defined from \eqref{Boundary Flow and Effort}-\eqref{Balance Equation}, and also, for all $x \in D(A(t))$,
	\begin{align}
			A(t)x &= \mathcal{J}_{l(t)} (\mathcal{Q}_{l(t)} x) = \begin{bmatrix}
				0 & \mathbf{d}_{l(t)} \\
				- \mathbf{d}_{l(t)}^{\ast} & 0 
			\end{bmatrix} (\mathcal{Q}_{l(t)} x).
		\label{Operator Family A(t)-II}
	\end{align} 	
 The operators $\mathbf{d}_{l(t)}$ and $\mathbf{d}_{l(t)}^{\ast}$ constituting the matrix differential operator $\mathcal{J}_{l(t)} \colon D(\mathcal{J}_{l(t)}) \subset X \to X$ given by
	\begin{align}
		\begin{split}
			D(\mathcal{J}_{l(t)}) &= D(\mathbf{d}_{l(t)}^{\ast}) \times D(\mathbf{d}_{l(t)}), \\
			\mathcal{J}_{l(t)}x &= \begin{bmatrix}
				0 & \mathbf{d}_{l(t)} \\
				- \mathbf{d}_{l(t)}^{\ast} & 0 
			\end{bmatrix} x, \quad x \in D(\mathcal{J}_{l(t)}), 
		\end{split}
		\label{Operator Jl(t)}
	\end{align} 
	are defined as in \eqref{Operator dl} and \eqref{Operator dl_ast} with respect to the time-varying color functions
	\begin{equation}
		c_{l}^-(z,t) = \begin{cases}
			1, & z \in [a,l(t)), \\
			0, & z \in [l(t),b],
		\end{cases}\quad t \geq 0,
		\label{Color Functions Depending on Moving Interface-I}
	\end{equation}
	\begin{equation}
		 c_{l}^+ (z,t) = \begin{cases} 
			0, & z \in [a,l(t)], \\
			1, & z \in (l(t),b],
		\end{cases}\quad t \geq 0.
		\label{Color Functions Depending on Moving Interface-II}
	\end{equation}
	Furthermore, the coercive operator $\mathcal{Q}_{l(t)} \in \mathcal{L}(X)$, $t \in [0, \tau]$, is given by 
	\begin{equation}
		\mathcal{Q}_{l(t)} = c_{l}^-(\cdot,t) \mathcal{Q}^- + c_{l}^+(\cdot,t) \mathcal{Q}^+.
		\label{Operator Ql(t)}
	\end{equation}
	We claim that, under certain conditions, $(A(t))_{t \in [0,\tau]}$ is a stable family on the state space $(X, \| \cdot \|_{\mathcal{Q}_0})$.
 Using \cite[Theorem 3.4]{KMMW23},
 and due to the equivalence  \eqref{Norm Equivalence} of the norms of the respective energy spaces $(X, \|\cdot\|_{\mathcal{Q}_{l(t)}})$, we can at least under natural conditions guarantee that the family $(A(t))_{t \in [0,\tau]}$ is a family of infinitesimal generators of uniformly bounded $C_0$-semigroups on $(X, \| \cdot \|_{\mathcal{Q}_0})$ (or any other state space $(X, \|\cdot\|_{\mathcal{Q}_{\hat{l}}})$ with $\hat{l} \in (a,b)$). However, as we will see, this does not prove yet the stability of this family. Thus, we need to impose more restrictive assumptions, namely:
	

	\begin{enumerate}[label =($\mathrm{A}$\arabic*)]
		\item \label{AssumptionA1} The matrix functions $\mathcal{Q}^{\pm} \in \mathcal{C}^1([a,b], \mathbb{R}^{2 \times 2})$ defining the coercive operator $\mathcal{Q}_{l(t)} \in \mathcal{L}(X)$ are diagonal, i.e., 
		\begin{align*}
		    \mathcal{Q}^{\pm}(z) = \begin{bmatrix}
		        q_{11}^{\pm}(z) & 0 \\
		        0 & q_{22}^{\pm}(z)
		    \end{bmatrix}, \quad z \in [a,b],
		\end{align*}
		and they satisfy 
      \begin{align*}
    \frac{q_{11}^+}{q_{11}^-}(0) = 1 \quad \text{and} \quad \frac{q_{11}^+}{q_{11}^-}(z) = \frac{q_{22}^+}{q_{22}^-}(z), \quad z \in [a,b]. 
    \end{align*}
		\item \label{AssumptionA2} $r =0 $, $\rank(W_B) = 2$, and $W_B \Sigma W_B^{\top} \geq 0$, where $\Sigma = \begin{bmatrix} 0 & I_2 \\
  I_2 & 0 
  \end{bmatrix} \in \R^{4 \times 4}$.
	\end{enumerate}
        Assumption \ref{AssumptionA1} states first that the Hamiltonian density function is a separable function, which is the case in numerous physical models such as the transmission line \cite{Gorrec2005}, the acoustic duct \cite{Kilian2022} or the isentropic gas model \cite{Diagne2013}.  Secondly, it states that the ratios of the coefficients defining the energy density functions of the subsystems are identical on the spatial domain $[a,b]$, and that they even coincide in $z = 0$. This happens for instance in the case of the interface separating two identical fluids, see Example 8 in \cite{Diagne2013}.
        
	Assumption \ref{AssumptionA2} requires first that $r =0 $, meaning that for any value of  $e_I(t)$, the variables $f_I(t) = 0$, implying that their power product is zero. This might happen, for instance, in the case of two gases separated by a mass-less piston which glides without dissipation. Hence, at the interface, power flow between the two subsystems can only occur by the motion of the interface. Also, the matrix inequality is precisely the condition that $A(t)$ generates a contraction semigroup on the respective energy space $(X, \|\cdot\|_{\mathcal{Q}_{l(t)}})$ for a \emph{frozen} interface position, see \cite[Theorem 5.14]{Kilian2022} \cite[Theorem 4.4]{KMMW23},
with $e_I(t) f_I(t)=0$ (as $r=0$) and the output conjugated to the velocity of the interface being the
\emph{discontinuity of energy density at the interface} \cite{Diagne2013}
\[
e_{l}=\left(-\mathcal{H}^{-}\left(l\right)+\mathcal{H}^{+}\left(l\right)\right).
\]

	To prove Kato-stability, we make use of an extension of a Lyapunov theorem \cite{Datko1968}.

	\begin{theorem}
		\label{Theorem Stability of A(t)}
		Let $0 < \tau < \infty$, and let $l \colon [0, \tau] \to (a,b)$ be continuously differentiable.
  Let the assumptions \ref{AssumptionA1}-\ref{AssumptionA2} hold. Consider the family $(A(t))_{t \in[0,\tau]}$ of port-Hamiltonian operators given by \eqref{Operator Family A(t)-I}-\eqref{Operator Ql(t)} defined on $(X, \|\cdot\|_{\mathcal{Q}_0})$. Then there exists some $\omega > 0$ such that for every $t \in [0, \tau]$ and for every $x \in D(A(t))$,
		\begin{equation}
			\langle A(t)x, x \rangle_{\mathcal{Q}_0} \leq \omega \| x\|_{\mathcal{Q}_0}^2.
			\label{Necessary Condition Stablity of Family}
		\end{equation} 
	\end{theorem}
	
\begin{IEEEproof}
	Let $t \in [0, \tau]$ be such that for the interface position we have $a < l(t) < 0$. In this case, the operator $\mathcal{J}_{l(t)}$ defined in \eqref{Operator Jl(t)} acts as the operator $P_1 \frac{d}{dz}$ on the subdomains $[a, l(t))$ and $(l(t),b]$, with $P_1$ defined in \eqref{Matrix P1}. As $A(t)$ generates a contraction semigroup on $(X, \|\cdot\|_{\mathcal{Q}_{l(t)}})$, it is dissipative, and for all $x \in D(A(t))$ it holds that
		\begin{align*}
			 \langle A(t)x,& x \rangle_{\mathcal{Q}_{l(t)}} = \frac{1}{2} \int_{a}^{b} x^{\top}(z) \mathcal{Q}_{l(t)}(z) [A(t) x](z) \, dz \\ 
			&= \frac{1}{2} \int_{a}^{l(t)} x^{\top}(z)  \mathcal{Q}^-(z) P_1 \frac{d}{dz} \big(\mathcal{Q}^-(z) x(z)\big) \, dz \\
			&\quad+ \frac{1}{2} \int_{l(t)}^{0} x^{\top}(z)  \mathcal{Q}^+(z) P_1 \frac{d}{dz} \big(\mathcal{Q}^+(z) x(z)\big) \, dz \\
			&\quad + \frac{1}{2} \int_{0}^{b} x^{\top}(z)  \mathcal{Q}^+(z) P_1 \frac{d}{dz} \big(\mathcal{Q}^+(z) x(z)\big) \, dz \leq 0.
		\end{align*}
	By dissipativity and since $\mathcal{Q}^{\pm} \in \mathcal{C}^1([a,b], \mathbb{R}^{2 \times 2})$, we have:
		\begin{align*}
			 \langle &A(t)x, x \rangle_{\mathcal{Q}_0} 
			= \langle A(t)x,x \rangle_{\mathcal{Q}_{l(t)}} \\
&+ \frac{1}{2} \int_{l(t)}^{0} x^{\top}(z) \big(\mathcal{Q}^-(z) - \mathcal{Q}^+(z)\big) P_1 \frac{d}{dz} (\mathcal{Q}^+x)(z) \, dz \\
			&\leq \frac{1}{2} \int_{l(t)}^{0} x^{\top}(z)\left[ \big(\mathcal{Q}^-(z) - \mathcal{Q}^+(z)\big)P_1 \frac{d}{dz}\mathcal{Q}^+(z)\right] x(z) \, dz \\
			& + \frac{1}{2} \int_{l(t)}^{0}  x^{\top}(z) \left[ \big(\mathcal{Q}^-(z) - \mathcal{Q}^+(z)\big)P_1 \mathcal{Q}^+(z) \right] \frac{d}{dz} x(z) \, dz.
		\end{align*}
	
	We may choose $\omega_1^- > 0$ large enough such that 
		\begin{equation*}
			\big(\mathcal{Q}^-(z) - \mathcal{Q}^+(z)\big)P_1 \frac{d}{dz}\mathcal{Q}^+(z) \leq \omega_1^- \mathcal{Q}^-(z), \quad z \in [a,b], 
		\end{equation*}
	i.e., there exists some $\omega_1^-$ such that for all $z \in [a,b]$, the symmetric part of 
	\begin{equation*}
		\omega_1^- \mathcal{Q}^-(z) - (\mathcal{Q}^-(z) - \mathcal{Q}^+(z))P_1 \frac{d}{dz}\mathcal{Q}^+(z)
	\end{equation*}
 	is positive semi-definite. Now, let us define the matrix function $\tilde{Q} \colon [a,b] \to \mathbb{R}^{2 \times 2}$ by
		\begin{align}
			\label{Theorem Stability of A(t) - Matrix TildeQ}
			\tilde{Q}(z) 
			&=  (\mathcal{Q}^-(z) - \mathcal{Q}^+(z))P_1 \mathcal{Q}^+(z) \nonumber\\
			&= - \begin{bmatrix}
				0 & q_{22}^+ (q_{11}^- - q_{11}^+) \\
				q_{11}^+(q_{22}^- - q_{22}^+)  & 0 
			\end{bmatrix}(z). 
		\end{align}
		By assumption \ref{AssumptionA1}, we have $\mathcal{Q}^+, \mathcal{Q}^- \in \mathcal{C}^1([a,b], \mathbb{R}^{2 \times 2}$), and so the matrix function $\tilde{Q}$ is continuously differentiable as well. Moreover, by \ref{AssumptionA1}, it is easy to check that $\tilde{Q}(z)$ is symmetric for all $z \in [a,b]$. Now integration by parts yields
		\begin{align*}
			& \quad \int_{l(t)}^{0} x^{\top}(z) \tilde{Q}(z) \frac{d}{dz} x(z) \, dz = \frac{1}{2} \big[x^{\top}(z)  \tilde{Q}(z) x(z) \big]_{l(t)^+}^{0} \\
   & - \frac{1}{2} \int_{l(t)}^{0} x^{\top}(z) \frac{d}{dz} \tilde{Q}(z) x(z) \, dz.
		\end{align*}
		Once again, we may choose $\omega_2^- > 0$ large enough such that
		\begin{equation*}
			- \frac{d}{dz} \tilde{Q}(z) \leq \omega_2^- \mathcal{Q}^-(z), \quad z \in [a,b]. 
		\end{equation*}
		Altogether, we obtain for all $x \in D(A(t))$, 
		\begin{align}	
				&\langle A(t)x, x \rangle_{\mathcal{Q}_0}	\leq  \frac{\omega_1^-}{2} \int_{l(t)}^{0} x^{\top}(z) \mathcal{Q}^-(z) x(z) \, dz \label{Theorem Auxiliary - A(t) is Stable - Proof Estimate 1} \\
				 &+ \frac{1}{4} \big[x^{\top}(z)  \tilde{Q}(z) x(z) \big]_{l(t)^+}^{0}
				+ \frac{\omega_2^-}{4} \int_{l(t)}^{0}  x^{\top}(z)\mathcal{Q}^-(z)  x(z) \, dz. \nonumber
		\end{align}     
		
		As $x \in D(A(t))$, the interface relation $f_I(t) = re_I(t)$ holds. Recalling the definition of the interface variables \eqref{Continuity Equation}-\eqref{Balance Equation}, invoking the assumption that $r = 0$ together with the assumption that $\mathcal{Q}^{\pm}$ are diagonal matrix functions, we obtain  
		\begin{align*}
            0 = f_I(t) 
						&= f_{I, \mathcal{Q}_{l(t)}x}
            = \left( \mathcal{Q}_{l(t)} x \right)_2 (l(t)^+)\\
            &= \left( \mathcal{Q}^+ x \right)_2 (l(t)^+)
            =  q_{22}^+(l(t)^+) x_2 (l(t)^+).
		\end{align*}
		As $q_{22}^+(l(t)^+)\neq 0$, we obtain that $x_2 (l(t)^+) = 0$. A similar computation shows that $x_2 (l(t)^-)=0$.
		
		Furthermore, by assumption \ref{AssumptionA1}, $\mathcal{Q}^{\pm}$ coincide in $z = 0$. 
		Consequently, the expression $\big[x^{\top}(z)  \tilde{Q}(z) x(z) \big]_{l(t)^+}^{0}$, i.e.,
	\begin{align*}
& x_1(l(t)^+) x_2(l(t)^+)  \left[ q_{22}^+ (q_{11}^- - q_{11}^+) + 	q_{11}^+(q_{22}^- - q_{22}^+) \right] (l(t)) \\
&-  x_1(0)x_2(0) \left[  q_{22}^+ (q_{11}^- - q_{11}^+) + 	q_{11}^+(q_{22}^- - q_{22}^+) \right] (0),
		\end{align*}
		vanishes. With the choice $\omega^- = \max\{\omega_1^-, \frac{\omega_2^-}{2}\}$ we deduce from \eqref{Theorem Auxiliary - A(t) is Stable - Proof Estimate 1} that for all $x \in D(A(t))$, where $t \in [0, \tau]$ is chosen such that $a < l(t) < 0$, it holds that
		\begin{equation*}
			\langle A(t)x, x \rangle_{\mathcal{Q}_0} \leq \omega^- \|x\|_{\mathcal{Q}_0}^2. 
		\end{equation*}
		Now, let $t \in [0, \tau]$ such that for the interface position we have $0 < l(t) < b$. Analogously to the previous case, we have for all $x \in D(A(t))$, 
		\begin{align*}
			& \quad \langle A(t)x, x \rangle_{\mathcal{Q}_0}\\
			&\leq \frac{1}{2} \int_{0}^{l(t)} x^{\top}(z) \left[ \big(\mathcal{Q}^+(z) - \mathcal{Q}^-(z)\big)P_1 \frac{d}{dz} \mathcal{Q}^-(z) \right] x(z) \, dz \\
			& + \frac{1}{2} \int_{0}^{l(t)} x^{\top}(z) \left[ \big(\mathcal{Q}^+(z) - \mathcal{Q}^-(z)\big)P_1 \mathcal{Q}^-(z)  \right] \frac{d}{dz} x(z) \, dz. 
		\end{align*}
		One can readily see that the previous arguments apply, and so there exists some $\omega^+ > 0$ such that 
		\begin{equation*}
			\langle A(t)x, x \rangle_{\mathcal{Q}_0} \leq \omega^+ \|x\|_{\mathcal{Q}_0}^2
		\end{equation*}
		for all $x \in D(A(t))$ with $0 < l(t) < b$. If $l(t) = 0$ for some $t \in [0,\tau]$, then $A(t)$ is dissipative on $(X, \|\cdot\|_{\mathcal{Q}_0})$. Altogether, the choice $\omega = \max\{\omega^-, \omega^+\}$ yields the estimate \eqref{Necessary Condition Stablity of Family} for all $t \in [0, \tau]$ and for all $x \in D(A(t))$. 
\end{IEEEproof}
	

    Using Theorem \ref{Theorem Stability of A(t)}, we have the following:
	\begin{corollary}
		\label{Corollary Stability of A(t)}
		Let $0 < \tau < \infty$, and let $l \colon [0, \tau] \to (a,b)$ be continuously differentiable. Let the assumptions \ref{AssumptionA1}-\ref{AssumptionA2} hold. Then each $A(t)$ generates a quasi-contractive semigroup, and the family $(A(t))_{t \in[0,\tau]}$ of port-Hamiltonian operators given by \eqref{Operator Family A(t)-I}-\eqref{Operator Ql(t)} is Kato-stable on $(X, \|\cdot\|_{\mathcal{Q}_0})$.
	\end{corollary}

\begin{IEEEproof}
 By \eqref{Necessary Condition Stablity of Family}, for all $x \in D(A(t))$ and $s \geq 0$
 \begin{align*}
			\frac{1}{2} \frac{d}{ds} \| S_t(s)x \|_{\mathcal{Q}_0}^2 \hspace{-0.055 cm}= \hspace{-0.055 cm} \langle A(t)S_t(s)x, S_t(s)x \rangle_{\mathcal{Q}_0} \hspace{-0.055 cm} \leq \omega \|S_t(s)x\|_{\mathcal{Q}_0}^2.
	\end{align*} 
Gr{\"o}nwall's Lemma (see, e.g., \cite[Section 2.4]{Teschl2012}) yields
	\begin{equation*}
	    \| S_t(s)x \|_{\mathcal{Q}_0}^2 \leq \| x\|_{\mathcal{Q}_0}^2 e^{2 \omega s}, \quad s \geq 0. 
	\end{equation*}
	As $D(A(t))$ is a dense subset of $X$, we conclude that
	\begin{equation*}
		\|S_t(s)\|_{\mathcal{L}(X, \|\cdot\|_{\mathcal{Q}_0})} \leq e^{\omega s}, \quad s \geq 0. 
	\end{equation*}
	Altogether, for all $t \in [0, \tau]$,  $(S_t(s))_{s \geq 0}$ is a quasicontractive semigroup on $(X, \|\cdot\|_{\mathcal{Q}_0})$. 
 
By the Hille-Yosida Theorem (see \cite[Theorem II.3.8]{Engel2000}) we have $(\omega, \infty) \subset \rho(A(t))$ for all $t \in [0, \tau]$. Furthermore, for every finite sequence $0 \leq t_1 \leq t_2 \leq \ldots \leq t_k \leq \tau$, $k \in \mathbb{N}$, and for all $\lambda > \omega$ it holds that
	\begin{align*}
		\left\| \prod_{j = 1}^{k} R(\lambda, A(t_j)) \right\|_{ \mathcal{L}(X, \|\cdot\|_{\mathcal{Q}_0})} &\leq \prod_{j=1}^{k} \left\| R(\lambda, A(t_j)) \right\|_{ \mathcal{L}(X, \|\cdot\|_{\mathcal{Q}_0})}\\
		&\leq \prod_{j=1}^{k} \frac{1}{\lambda - \omega} = \frac{1}{(\lambda - \omega)^k}.
	\end{align*}
 Thus, the family $(A(t))_{t \in [0, \tau]}$ is Kato-stable with stability constants $M=1$ and $\omega > 0$.
\end{IEEEproof}

	Note that, while the stability of the family $(A(t))_{t \in [0, \tau]}$ is a crucial result, it is not sufficient for well-posedness of the corresponding evolution problem \eqref{Evolution Problem}, see, e.g., \cite[Chapter 5]{Pazy1983}, \cite{Kobayasi1979}, \cite{Nickel1998}. In fact, another important prerequisite is the existence of a densely and continuously embedded Banach space $Y \subset X$ satisfying
	\begin{equation}
	    Y \subset \bigcap_{t \in [0, \tau]} D(A(t)).
	\end{equation}
    However, due to the fact that the interface conditions \eqref{Continuity Equation}-\eqref{Balance Equation} and $f_{I, \mathcal{Q}_{l} x} = r e_{I, \mathcal{Q}_{l} x}$ (with $r = 0$) have to be satisfied for all interface positions in $\range(l)$, such a subspace only exists for the trivial case that $l(t)$ remains constant over time. Hence, standard tools for verifying that the family $(A(t))_{t \in [0, \tau]}$ generates an evolution family fail immediately.
    
	\section{A counterexample}
 \label{sec:counterexample}
In this section, we show that in the absence of condition \ref{AssumptionA1}, we cannot expect that the strategy of proof for Kato-stability works. While this does not show that the family is indeed not Kato-stable, it indicates obstacles for proving this.

\begin{example}
Here we show that in the absence of assumption \ref{AssumptionA1} the result of Theorem~\ref{Theorem Stability of A(t)} does not hold. Namely, 
there exist $\cQ^{\pm}\in \mathcal{C}^1\bigl(\left[-2,2\right],\R^{2\times 2} \bigr)$ both diagonal, coercive and bounded and such that for the induced operator $\cQ_{-1}$ there exists $x_k \in \mathcal{C}^{\infty}\bigl(\left[-\frac{3}{4}, -\frac{1}{4}\right], \R^2\bigr),k\in \N$, with $\|x_k\|_2\leq 1, k\in \N$, and
\begin{equation*}
\langle A_{\cQ_{-1}}x_k, x_k\rangle_{\cQ_{0}}\rightarrow\infty\,.
\end{equation*}
Note that as the support of $x_k$ is contained in $\left[-\frac{3}{4}, -\frac{1}{4}\right]$, the boundary conditions at $z=-2, -1,0,2$ are automatically satisfied so this does not play a role in the sequel.

For any smooth $x\in X$ with $\text{supp}\,x\subset \left[-\frac{3}{4}, -\frac{1}{4}\right]$, we have $x\in D\bigl(A_{\cQ_{-1}} \bigr)$ and 
\begin{multline*}
\langle A_{\cQ_{-1}}x, x\rangle_{\cQ_0}  =
\int_{-\frac{3}{4}}^{-\frac{1}{4}}x^{\top} \cQ^-P_1\frac{d}{dz}\bigl(\cQ^+x\bigr) \, dz\\
=\int_{-\frac{3}{4}}^{-\frac{1}{4}}x^{\top}
 \cQ^-P_1\cQ^{+'}x \, dz+\int_{-\frac{3}{4}}^{-\frac{1}{4}}x^{\top}\cQ^- P_1\cQ^+x' \, dz\,.
\end{multline*}

The first summand can be bounded in absolute value in terms of $\|x\|_2$. We thus only need to consider the second term.
%
%
For simplicity, write $\cQ^-=\text{diag}\,\bigl(\tilde{q}_1, \tilde{q}_2\bigr)$ and $\cQ^+=\text{diag}\,\bigr(\tilde{r}_1, \tilde{r}_2\bigl)$. Then the second term becomes
\begin{align*}
& =\int_{-\frac{3}{4}}^{-\frac{1}{4}}x^{\top}\begin{bmatrix}
0 & -\tilde{q}_1\tilde{r}_2\\
-\tilde{q}_2\tilde{r}_1 & 0
\end{bmatrix}
x' \, dz\\
 & =\int_{-\frac{3}{4}}^{-\frac{1}{4}}-\tilde{q}_2\tilde{r}_1 x_2(z)x_1'(z)-\tilde{q}_1\tilde{r}_2x_1(z)x_{2}'(z) \, dz\,.
\end{align*}

We now assume that for some $\varepsilon >0$, $\tilde{q}_1\tilde{r}_2 =\varepsilon$ and
$\tilde{q}_2\tilde{r}_1 =\varepsilon+\chi_{[\xi_1,\xi_2]}+\eta_1$,
where $[\xi_1, \xi_2]\subset \bigl(-\frac{3}{4},-\frac{1}{2}
\bigr)$, $\chi_{[\xi_1,\xi_2]}$ is the characteristic function of $[\xi_1, \xi_2]$, 
and
$\eta_1$ is  a function guaranteeing smoothness of
$\tilde{q}_2\tilde{r}_{1}$ and with small $L^2$-norm.
The constant integral (from $\varepsilon$) evaluates to zero, so 
\begin{equation*}
=-\int_{\xi_1}^{\xi_2}x_2(z)x_{1}'(z) \, dz-\int_{\xi_3}^{\xi_4}x_1(z)x_{2}'(z) \, dz+\delta(\eta_1)\,.
\end{equation*}

Now choose $x_2$ smooth with 
$x_2(z) =-1, z\in [\xi_1,\xi_2]$
and $\|x_2\|_2=\frac{1}{3}$. Then we can continue
\begin{equation*}
=\int_{\xi_1}^{\xi_2}x_{1}'(z) \, dz+\delta=x_1(\xi_2)-x_1(\xi_1)+\delta\,.
\end{equation*}
It is well known that the last expression can be arbitrarily large for smooth functions $x_1$ with bounded $L_2$-norm.
\end{example}

 \section{Conclusions}
 \label{sec:conclusions}

 We have considered port-Hamiltonian systems on a one-dimensional domain with a moving interface. For a particular case, it can be shown that the family of generators corresponding to fixed interfaces is Kato-stable. Even in the  setting of this paper, the question of generation of an evolution family remains open. An example shows that the argument for Kato-stability fails in more general situations, even though the associated generators for fixed interfaces are benign and can be treated by the methods developed in \cite{Kilian2022}.

  For future elaboration, one could consider the additional port variables associated with the velocity of the interface as well as general port variables at the interface defined as in (\ref{Boundary Flow and Effort}) and consider some passivity properties of interconnection models at the interface.

In \cite{KMMW23}, using semigroup methods, we have fully characterized boundary conditions making a port-Hamiltonian system with an interface a well-posed system governed by a contraction semigroup.
Alternatively, the well-posedness analysis could be approached using the methods of the theory of PDEs with moving boundaries, which is well-developed for the case of parabolic equations \cite{Cra84, PrS16}, but is also analyzed for hyperbolic systems \cite{WMW21}.

\bibliographystyle{IEEEtran}
\bibliography{references}

\begin{thebibliography}{10}
\providecommand{\url}[1]{#1}
\csname url@samestyle\endcsname
\providecommand{\newblock}{\relax}
\providecommand{\bibinfo}[2]{#2}
\providecommand{\BIBentrySTDinterwordspacing}{\spaceskip=0pt\relax}
\providecommand{\BIBentryALTinterwordstretchfactor}{4}
\providecommand{\BIBentryALTinterwordspacing}{\spaceskip=\fontdimen2\font plus
\BIBentryALTinterwordstretchfactor\fontdimen3\font minus
  \fontdimen4\font\relax}
\providecommand{\BIBforeignlanguage}[2]{{%
\expandafter\ifx\csname l@#1\endcsname\relax
\typeout{** WARNING: IEEEtran.bst: No hyphenation pattern has been}%
\typeout{** loaded for the language `#1'. Using the pattern for}%
\typeout{** the default language instead.}%
\else
\language=\csname l@#1\endcsname
\fi
#2}}
\providecommand{\BIBdecl}{\relax}
\BIBdecl

\bibitem{Duindam2009}
V.~Duindam, A.~Macchelli, S.~Stramigioli, and H.~Bruyninckx, Eds.,
  \emph{Modeling and {C}ontrol of {C}omplex {P}hysical {S}ystems}, 2009.

\bibitem{Jacob2012}
B.~Jacob and H.~J. Zwart, \emph{Linear {P}ort-{H}amiltonian {S}ystems on
  {I}nfinite-dimensional {S}paces}.\hskip 1em plus 0.5em minus 0.4em\relax
  Birkhäuser/Springer Basel AG, Basel, 2012.

\bibitem{Schaft2014}
A.~J. van~der Schaft and D.~Jeltsema, ``Port-{H}amiltonian systems theory: an
  introductory overview,'' \emph{Foundations and Trends® in Systems and
  Control}, vol.~1, no. 2-3, pp. 173--378, 2014.

\bibitem{Arnold1978}
V.~I. Arnol'd, \emph{Mathematical {M}ethods of {C}lassical {M}echanics}.\hskip
  1em plus 0.5em minus 0.4em\relax Springer, New York, 1978.

\bibitem{Curtain2020}
R.~F. Curtain and H.~J. Zwart, \emph{Introduction to {I}nfinite-{D}imensional
  {L}inear {S}ystems {T}heory}.\hskip 1em plus 0.5em minus 0.4em\relax
  Springer, New York, 2020.

\bibitem{Gorrec2005}
Y.~L. Gorrec, H.~J. Zwart, and B.~M. Maschke, ``Dirac structures and boundary
  control systems associated with skew-symmetric differential operators,''
  \emph{SIAM Journal on Control and Optimization}, vol.~44, no.~5, pp.
  1864--1892, 2005.

\bibitem{KMMW23}
A.~Kilian, A.~Mironchenko, B.~Maschke, and F.~Wirth, ``Infinite-dimensional
  port-{H}amiltonian systems with a stationary interface,'' \emph{Submitted to
  SIAM J. Control Optim.}, 2023, arXiv:2301.08967.

\bibitem{Kilian2022}
A.~Kilian, ``Infinite-{D}imensional {P}ort-{H}amiltonian {S}ystems with a
  {M}oving {I}nterface,'' Master's thesis, Chair of Dynamical Systems,
  University of Passau, 2022, arXiv:2301.07344.

\bibitem{Diagne2013}
M.~Diagne and B.~M. Maschke, ``Port {H}amiltonian formulation of a system of
  two conservation laws with a moving interface,'' \emph{European Journal on
  Control}, vol.~19, no.~6, pp. 495--504, 2013.

\bibitem{Pazy1983}
A.~Pazy, \emph{Semigroups of {L}inear {O}perators and {A}pplications to
  {P}artial {D}ifferential {E}quations}.\hskip 1em plus 0.5em minus 0.4em\relax
  Springer, New York, 1983.

\bibitem{Engel2000}
K.-J. Engel and R.~Nagel, \emph{One-{P}arameter {S}emigroups for {L}inear
  {E}volution {E}quations}.\hskip 1em plus 0.5em minus 0.4em\relax Springer,
  New York, 2000.

\bibitem{Datko1968}
R.~Datko, ``An extension of a theorem of {A}. {M}. {L}yapunov to semigroup
  operators,'' \emph{J. Math. Analysis Appl.}, vol.~24, pp. 290--295, 1968.

\bibitem{Teschl2012}
G.~Teschl, \emph{Ordinary {D}ifferential {E}quations and {D}ynamical
  {S}ystems}.\hskip 1em plus 0.5em minus 0.4em\relax American Mathematical
  Society, Providence, 2012.

\bibitem{Kobayasi1979}
K.~Kobayasi, ``On a theorem for linear evolution equations of hyperbolic
  type,'' \emph{J. Math. Soc. Japan}, vol.~31, no.~4, pp. 647--654, 1979.

\bibitem{Nickel1998}
G.~Nickel and R.~Schnaubelt, ``An extension of {K}ato's stability condition for
  nonautonomous {C}auchy problems,'' \emph{Taiwanese Journal of Mathematics},
  vol.~2, no.~4, pp. 483--496, 1998.

\bibitem{Cra84}
J.~Crank, \emph{Free and moving boundary problems}.\hskip 1em plus 0.5em minus
  0.4em\relax Oxford University Press, USA, 1984.

\bibitem{PrS16}
J.~Pr{\"u}ss and G.~Simonett, \emph{Moving interfaces and quasilinear parabolic
  evolution equations}.\hskip 1em plus 0.5em minus 0.4em\relax Springer, 2016,
  vol. 105.

\bibitem{WMW21}
J.~Wurm, L.~Mayer, and F.~Woittennek, ``Feedback control of water waves in a
  tube with moving boundary,'' \emph{European Journal of Control}, vol.~62, pp.
  151--157, 2021.

\end{thebibliography}

\end{document}